\documentclass[11pt,a4paper]{amsart}

\usepackage[latin2]{inputenc}
\usepackage{t1enc}
\usepackage{amsmath, amsthm, amssymb}
\usepackage{graphicx}
\usepackage{epstopdf}

\usepackage{color}

\newtheorem{Thm}{Theorem}[section]
\newtheorem{Prop}[Thm]{Proposition}
\newtheorem{Lem}[Thm]{Lemma}
\newtheorem{Cor}[Thm]{Corollary}
\newtheorem{Res}[Thm]{Result}

\theoremstyle{definition}

\newtheorem{Ex}[Thm]{Example}
\newtheorem{Rem}[Thm]{Remark}

\numberwithin{equation}{section}

\newenvironment{Proof}{\rm \trivlist\item[\hskip \labelsep{\bf
Proof.\quad}]}{\hfill\qed\par\medskip\endtrivlist}

\newcommand{\slU}{\mathbf{U}}

\newcommand{\slAb}{\mathbf{Ab}}
\newcommand{\slC}{\mathbf{C}}
\newcommand{\slD}{\mathbf{D}}
\newcommand{\slI}{\mathcal{I}}
\newcommand{\slX}{\mathcal{X}}
\newcommand{\slY}{\mathcal{Y}}
\newcommand{\ol}{\overline}
\newcommand{\ul}{\underline}

\newcommand{\rmE}{\mathrm{E}}
\newcommand{\rmX}{\mathrm{X}}
\newcommand{\rmY}{\mathrm{Y}}
\newcommand{\rmZ}{\mathrm{Z}}

\newcommand{\xc}{\slX_\mathrm{cl}(\Gamma^I)} 
\newcommand{\maxM}{\max M^-} 

\def\Sub{\operatorname{Sub}}
\def\ClSub{\operatorname{ClSub}}
\def\ClSubfc{\operatorname{ClSub_{\mathrm fc}}}

\def\Ob{\operatorname{Ob}}
\def\Arr{\operatorname{Arr}}

\def\act#1#2{{^{#1}\kern -2pt {#2}}}
\def\acth#1#2{{^{#1}\kern -1pt {#2}}}

\def\cat#1{F_{g#1}(\Gamma)} 
\def\cl#1{{#1}^\mathrm{cl}} 
\def\Cc#1{C_{#1}^\mathrm{cl}} 
\def\Pc#1{P_{#1}^\mathrm{cl}} 
\def\sp#1{\langle{#1}\rangle} 
\def\spc#1{\langle{#1}\rangle^\mathrm{cl}} 
\def\lbl{\lambda}
\def\ta{\tilde{a}}

\def\te{\tilde{e}}
\def\tf{\tilde{f}}

\def\tu{\tilde{u}}
\def\tA{\tilde{A}}

\begin{document}

\title[$F$-inverse covers of inverse monoids]{On $F$-inverse covers of finite-above inverse monoids}

\author{N\'ora Szak\'acs}
\address{Bolyai Institute, University of Szeged, Aradi v\'ertan\'uk tere 1, Szeged, Hungary, H-6720; fax: +36 62 544548}
\email{szakacsn@math.u-szeged.hu}

\author{M\'aria B.\ Szendrei}
\address{Bolyai Institute, University of Szeged, Aradi v\'ertan\'uk tere 1, Szeged, Hungary, H-6720; fax: +36 62 544548}
\email{m.szendrei@math.u-szeged.hu}

\date{Dec 31, 2015}

\thanks{Research supported by the Hungarian National Foundation for Scientific 
Research grant no.~K83219, K104251,
and by the European Union, under the project 
no.\ T\'{A}MOP-4.2.2.A-11/1/KONV-2012-0073 and T\'{A}MOP-4.2.2.B-15/1/KONV-2015-0006.
}

\begin{abstract}
Finite-above inverse monoids are a common generalization of finite inverse monoids and Margolis--Meakin expansions of groups.
Given a finite-above $E$-unitary inverse monoid $M$ and a group variety $\slU$, 
we find a condition for $M$ and $\slU$, involving a construction of descending chains of graphs, which is equivalent to
$M$ having an $F$-inverse cover via $\slU$.
In the special case where $\slU=\slAb$, the variety of Abelian groups, we apply 
this condition
to get a simple sufficient condition for $M$ to 
have no
$F$-inverse cover via $\slAb$, 
formulated by means of the natural partial order and the least group congruence of $M$.
\end{abstract}

\maketitle

\section{Introduction}

An {\em inverse monoid} is a monoid $M$ with the
property that for each $a \in M$ there exists a unique element
$a^{-1} \in M$ (the inverse of $a$) such that $a =aa^{-1}a$ and
$a^{-1} = a^{-1}aa^{-1}$. The {\em natural partial order} on an inverse monoid $M$ is defined
as follows: $a \leq b$ if $a = eb$ for some idempotent $e \in M$.  
Each inverse monoid admits a smallest group congruence which is usually denoted by $\sigma$.
Inverse monoids appear in various areas of mathematics due to their role in the description of partial 
symmetries (see \cite{lo} for this approach).

An inverse monoid is called $F$-inverse if each class of the least group congruence has a greatest element with respect to the natural partial order.
For example, free inverse monoids are $F$-inverse due to the fact that the Cayley graph of a free group is a tree.
This implies that every inverse monoid has an $F$-inverse cover, i.e., every inverse monoid is an 
(idempotent separating) homomorphic image of an $F$-inverse monoid.
This is the only known way to produce an $F$-inverse cover for any inverse monoid, but it constructs an infinite cover even for
finite inverse monoids.
It is natural to ask whether each finite inverse monoid has a finite $F$-inverse cover.
This question was first formulated by K.\ Henckell and J.\ Rhodes \cite{HR}, when they observed that 
an affirmative answer would imply an affirmative answer for the pointlike conjecture for inverse monoids. 
The latter conjecture was proved by C.\ J.\ Ash \cite{A}, but the finite $F$-inverse cover problem is still open.

K.\ Auinger and the second author reformulate the finite $F$-inverse cover problem in \cite{ASzM} as follows.
First they reduce the question of the existence of a finite $F$-inverse cover for all finite inverse monoids 
to the existence of a generator-preserving dual premorphism 
--- a map more general than a homomorphism ---
from a finite group to a Margolis--Meakin expansion of a finite group. 
The Margolis--Meakin expansion $M(G)$ of a group $G$ is obtained from a Cayley graph of 
$G$ similarly to how a free inverse monoid is obtained from the Cayley graph of a free group. 
If there is a dual premorphism $H \to M(G)$ from a finite group $H$, 
then $H$ is an extension of some group $N$ by $G$, and, for any group variety $\slU$ containing $N$,
there is a dual premorphism from the `most general' extension of 
a member of $\slU$ by $G$. 
This group, denoted by $G^\slU$, can also be obtained from the Cayley graph of $G$ and the variety $\slU$. 
It turns out that in fact the existence of a dual premorphism $G^\slU \to M(G)$ only depends on the Cayley graph of $G$ and the 
variety $\slU$. 
This leads in \cite{ASzM} to an equivalent of the finite $F$-inverse cover problem formulated by means of graphs and 
group varieties.

The motivation for our research was to obtain similar results for a class of inverse monoids significantly larger than Margolis-Meakin expansions.
Studying the dual premorphisms $(M/\sigma)^\slU \to M$, where $M$ is an inverse monoid and $\slU$ is a group variety,
we generalize the main construction of \cite{ASzM} for a class called finite-above $E$-unitary inverse monoids, which contains
all finite $E$-unitary inverse monoids and the Margolis--Meakin expansions of all groups. 
This leads to a condition
--- involving, as in \cite{ASzM}, 
a process of constructing descending chains
of graphs --- on a finite-above $E$-unitary inverse monoid $M$ and a group variety $\slU$ 
that is equivalent to $M$ having an $F$-inverse cover over $(M/\sigma)^\slU$,
or, as we shall briefly say later on, an $F$-inverse cover via $\slU$.
Our condition restricts to the one in \cite{ASzM} if $M$ is a Margolis--Meakin expansion of a group.
As an illustration, we apply this process to find a sufficient condition on the natural partial order and the least group congruence $\sigma$ of $M$ 
for $M$ to have no $F$-inverse cover 
via $\slAb$, the variety of Abelian groups.

In Section 2, we introduce some of the structures and definitions needed later, particularly those from the theory of inverse categories 
acted upon by groups. 
These play an important role in our approach, especially the derived categories of the natural morphisms of $E$-unitary inverse monoids onto
their greatest group images and some categories arising from Cayley graphs of groups. 
In the first part of Section 3, we introduce the inverse monoids which play the role of the 
Margolis--Meakin expansions of groups in our paper.
The second part of Section 3 contains the main result of the paper (Theorem \ref{main}). 
Section 4 is devoted to giving a sufficient condition in Theorem \ref{Abelmain} for a
finite-above $E$-unitary inverse monoid to have no $F$-inverse cover via the variety of Abelian groups.

\section{Preliminaries}

In this section we recall the notions and results needed in the paper.
For the undefined notions and notation, the reader is referred to \cite{lo} and \cite{pet}.

\subsection*{$A$-generated inverse monoids}
Let $M$ be an inverse monoid (in particular, a group) and $A$ an arbitrary set.
We say that $M$ is an {\em $A$-generated inverse monoid} ({\em $A$-generated group}) if a map $\epsilon_M\colon A\to M$
is given such that $A\epsilon_M$ generates $M$ as an inverse monoid (as a group).
If $\epsilon_M$ is injective, then we might assume that $A$ is a subset in $M$, i.e., $\epsilon_M$ is the inclusion map 
$A\to M$.

Let $M$ be an $A$-generated inverse monoid. Consider a set $A'$ disjoint from $A$ together with a bijection $'\colon A\to A'$.
Put $\ol A=A\cup A'$, and denote the free monoid on $\ol A$ by $\ol A^*$.
Then there is a unique homomorphism $\varphi\colon \ol A^*\to M$ such that $a\varphi=a\epsilon_M$ and $a'\varphi=(a\epsilon_M)^{-1}$
for every $a\in A$, and $\varphi$ is clearly surjective.
For any word $w\in \ol A^*$, we denote $w\varphi$ by $[w]_M$.
In particular, if $M$ is the relatively free group on $A$ in a given group variety $\slU$, then we write $[w]_{\slU}$ for $[w]_M$.
Recall that, for every $w,w_1\in \ol A^*$, we have $[w]_{\slU}=[w_1]_{\slU}$ if and only if the identity $w=w_1$
is satisfied in $\slU$.

\subsection*{Graphs and categories}
Throughout the paper, unless otherwise stated, by a {\em graph} we mean a directed graph.
Given a graph $\Delta$, its set of vertices and set of edges are denoted by $V_\Delta$ and $E_\Delta$ respectively.
If $e\in E_\Delta$, then $\iota e$ and $\tau e$ are used to denote the initial and terminal vertices of $e$, and if
$\iota e=i$, $\tau e=j$, then $e$ is called an $(i,j)$-edge.
The set of all $(i,j)$-edges is denoted by $\Delta(i,j)$, and for our later convenience, we put
$$\Delta(i,-)=\bigcup_{j\in V_\Delta}\Delta(i,j).$$
We say that $\Delta$ is {\em connected} if its underlying undirected graph is connected.
If a set $A$ and a map $E_\Delta\to A$ is given, then $\Delta$ is said to be {\em labelled by $A$}.
For example, the Cayley graph of an $A$-generated group $G$ is connected and labelled by $A$.
Moreover, if $1\notin A\epsilon_G$, then there are no loops in the Cayley graph.
A sequence $p=e_1e_2\cdots e_n\ (n\ge 1)$ of consecutive edges $e_1,e_2,\ldots,e_n$
(i.e., where $\tau e_i=\iota e_{i+1}\ (i=1,2,\ldots,n-1)$)
is called a {\em path} on $\Delta$ or, more precisely, an $(i,j)$-path if $i=\iota e_1$ and $j=\tau e_n$.
In particular, if $i=j$ then $p$ is also said to be a {\em cycle} or, more precisely, an {\em $i$-cycle}.
Moreover, for any vertex $i\in V_\Delta$, we consider an {\em empty $(i,i)$-path} ({\em$i$-cycle}) denoted by $1_i$.
A non-empty path (cycle) $p=e_1e_2\cdots e_n$ is called {\em simple} if the vertices $\iota e_1,\iota e_2,\ldots, \iota e_n$ 
are pairwise distinct and $\tau e_n\notin \{\iota e_2,\ldots, \iota e_n\}$.

As it is usual with paths on Cayley graphs of groups, we would like to allow traversing edges in the reverse direction.
More formally, we add a formal reverse of each edge to the graph, and consider paths on this extended graph as follows.
Given a graph $\Delta$, consider a set $E'$ disjoint from $E_\Delta$ together with a bijection
$'\colon E_\Delta\to E'$, and consider a graph $\Delta'$ where $V_{\Delta'}=V_\Delta$ and $E_{\Delta'}=E'$
such that $\iota e'=\tau e$ and $\tau e'=\iota e$ for every $e\in E_\Delta$.
Define $\ol\Delta$ to be the graph with $V_{\ol\Delta}=V_\Delta$ and $E_{\ol\Delta}=E_\Delta\cup E_{\Delta'}$.
Choosing the set $E_\Delta'$ to be $E_{\Delta'}$,
the paths on $\ol\Delta$ become words in $\ol{E_\Delta}^*$ where $\ol{E_\Delta}=E_\Delta\cup E_\Delta'$.
Most of our graphs in this paper have edges of the form $(i,a,j)$, where $i$ is the initial 
vertex, $j$ the terminal vertex, and $a$ is the label of the edge.
For such a graph $\Delta$, we choose $\ol\Delta$ as follows:
we consider a set $A'$ disjoint from $A$ together with a bijection $'\colon A\to A'$ (see the previous subsection), and we choose $\Delta'$ so that
$(i,a,j)'=(j,a',i)$ for any edge $(i,a,j)$ in $\Delta$.
Then $\ol\Delta$ is labelled by $\ol A$, and, 
given a (possibly empty) path $p=e_1e_2\cdots e_n$ on $\ol\Delta$, the labels of the edges $e_1,e_2,\ldots,e_n$
determine a word in $\ol A^*$.

Now we extend the bijection $'$ to paths in a natural way.
First, for every edge $f\in E_{\Delta'}$, define $f'=e$ where $e$ is the unique edge in $\Delta$ such that $e'=f$.
Second, put $1_i'=1_i\ (i\in V_\Delta)$ and, for every non-empty path $p=e_1e_2\cdots e_n$ on $\ol\Delta$,
put $p'=e'_ne'_{n-1}\cdots e'_1$.
If $p=e_1e_2\cdots e_n$ is a non-empty path on $\ol\Delta$, then {\em the subgraph $\sp p$ of $\Delta$ spanned by $p$}
is the subgraph consisting of all
vertices and edges $p$ traverses in either direction.
Obviously, we have $\sp{p'}=\sp p$ for any path $p$ on $\ol\Delta$.
The subgraph spanned by the empty path $1_i$ (consisting of the single vertex $i$) is
denoted by $\emptyset_i$, that is, $\sp{1_i}=\emptyset_i$.

Let $\Delta$ be a graph, and suppose that a partial multiplication is given on $E_\Delta$ in a way that,
for any $e,f\in E_\Delta$, the product $ef$ is defined if and only if $e$ and $f$ are consecutive edges.
If this multiplication is associative in the sense that $(ef)g=e(fg)$ whenever $e,f,g$ are consecutive, and
for every $i\in V_\Delta$, there exists (and if exists, then unique) edge $1_i$ with the property that
$1_ie=e,\ f1_i=f$ for every $e,f\in E_\Delta$ with $\iota e=i=\tau f$, then $\Delta$ is called a {\em (small) category}.
Later on, we denote categories in calligraphics.
For categories, the usual terminology and notation is different from those for graphs: 
instead of `vertex' and `edge', we use the terms `object' and `arrow', respectively, and if $\slX$ is a category, then, instead of 
$V_\slX$ and $E_\slX$, we write $\Ob\slX$ and $\Arr\slX$, respectively. 
Clearly, each monoid can be considered a one-object category. 
Therefore, later on, certain definitions and results formulated only for categories will be applied also for monoids.
Given a graph $\Delta$, we can easily define a category $\Delta^*$ as follows: 
let $\Ob\Delta^*=V_\Delta$, let $\Delta^*(i,j)\ (i,j\in \Ob\Delta^*)$ be the set of all $(i,j)$-paths on $\Delta$, 
and define the product of consecutive paths by concatenation.
The identity arrows will be the empty paths.
In the one-object case, this is just the usual construction of a free monoid on a set. 
In general, $\Delta^*$ has a similar universal property among categories, that is, it is the {\em free category on $\Delta$}.

A category $\slX$ is called a {\em groupoid} if, for each arrow $e\in\slX(i,j)$, there exists an arrow $f\in\slX(j,i)$ such that
$ef=1_i$ and $fe=1_j$.
Obviously, the one-object groupoids are just the groups, and, as it is well known for groups, the arrow $f$ is uniquely determined,
it is called the inverse of $e$ and is denoted $e^{-1}$.
By an {\em inverse category}, we mean a category $\slX$ where, for every arrow $e\in\slX(i,j)$, there exists a unique arrow $f\in\slX(j,i)$ such that
$efe=e$ and $fef=f$. This unique $f$ is also called the {\em inverse of $e$} and is denoted $e^{-1}$.
Clearly, each groupoid is an inverse category, and this notation does not cause confusion.
Furthermore, the one-object inverse categories are just the inverse monoids.
More generally, if $\slX$ is an inverse category (in particular, a groupoid), then $\slX(i,i)$ is an inverse monoid (a group) for every 
object $i$.
An inverse category $\slX$ is said to be {\em locally a semilattice} if $\slX(i,i)$ is a semilattice for every object $i$.
Similarly, given a group variety $\slU$, we say that $\slX$ is {\em locally in $\slU$} if $\slX(i,i)\in\slU$ for every object $i$. 
For an inverse category $\slX$ and a graph $\Delta$, if $\epsilon_\slX\colon \Delta\to\slX$ is a graph morphism,
then there is a unique category morphism $\varphi\colon \ol\Delta^*\to \slX$ such that 
$e\varphi=e\epsilon_\slX$ and $e'\varphi=(e\epsilon_\slX)^{-1}$
for every $e\in \Arr\slX$.
We say that $\slX$ is {\em $\Delta$-generated} if $\varphi$ is surjective.

The basic notions and properties known for inverse monoids have their analogues for inverse categories.
Given a category $\slX$, consider {\em the subgraph $\rmE(\slX)$ of idempotents}, where $V_{\rmE(\slX)}=\Ob{\slX}$ and
$E_{\rmE(\slX)}=\{h\in\Arr{\slX}:hh=h\}$.
Obviously, $E_{\rmE(\slX)}\subseteq\bigcup_{i\in\Ob{\slX}}\slX(i,i)$.
A category $\slX$ is an inverse category if and only if $\rmE(\slX)(i,i)$ is a semilattice for every object $i$, and,
for each arrow $e\in\slX(i,j)$, there exists an arrow $f\in\slX(j,i)$ such that $efe=e$ and $fef=f$.
Thus, given an inverse category $\slX$, $\rmE(\slX)$ is a subcategory of $\slX$, and we define a relation $\le$ on $\slX$ as follows: 
for any $e,f\in \Arr\slX$, let $e\le f$ if $e=fh$ for some $h\in\Arr{\rmE(\slX)}$.
The relation $\le $ is a partial order on $\Arr\slX$ called the {\em natural partial order} on $\slX$, and
it is compatible with multiplication.
Note that the natural partial order is trivial if and only if $\slX$ is a groupoid.

\subsection*{Categories acted upon by groups}
Now we recall several notions and facts from \cite{mar-pin}.

Let $G$ be a group and $\Delta$ a graph. We say that {\em $G$ acts on $\Delta$ (on the left)} if, 
for every $g\in G$, and for every vertex $i$ and edge $e$ in $\Delta$, 
a vertex $\acth gi$ and an edge $\acth ge$ is given such that the following are satisfied for any $g,h\in G$ and any $i\in V_\Delta$, $e\in E_\Delta$:
$$\acth 1i=i,\quad \acth h{(\acth gi)}=\acth{hg}i,\quad \acth 1e=e,\quad \acth h{(\acth ge)}=\acth{hg}e,$$
$$\iota\acth ge=\acth g{\iota e},\quad \tau\acth ge=\acth g{\tau e}.$$
An action of $G$ on $\Delta$ induces an action on the paths and an action on the subgraphs of $\Delta$ in a natural way:
if $g\in G$, $i\in V_\Delta$ and
$p=e_1e_2\cdots e_n$ is a non-empty path, then we put
$$\acth gp=\acth g{e_1}\acth g{e_2}\cdots \acth g{e_n},$$
and for an empty path, let $\acth g{1_i}=1_{\acth gi}$.
For any subgraph $\rmX$ of $\Delta$, define
$\acth g{\rmX}$ to be the subgraph whose sets of vertices and edges are 
$\{\acth gi: i\in V_\rmX\}$ and $\{\acth ge: e\in E_\rmX\}$ respectively, in particular, $\acth g{\emptyset_i}=\emptyset_{\acth gi}$.
The action of $G$ on $\Delta$ can be extended to $\ol{\Delta}$ also in a natural way by setting
$\acth g{e'}=(\acth ge)'$ for every $e\in E_\Delta$.
It is easy to check that the equality $\sp{\act gp}=\act g{\sp p}$ holds for every path $p$ on $\ol{\Delta}$.

By an {\em action of a group on a category $\slX$} we mean an action of $G$ on the graph $\slX$ which has the following additional properties:
for any object $i$ and any pair of consecutive arrows $e,f$, we have
$$\acth g{1_i}=1_{\acth gi},\quad \acth g{(ef)}=\acth ge\cdot \acth gf.$$
In particular, if $\slX$ is a one-object category, that is, a monoid, then this defines an action of a group on a monoid. 
We also mention that if $\Delta$ is a graph acted upon by a group $G$, then the induced action on the paths
defines an action of $G$ on the free categories $\Delta^*$ and $\ol{\Delta}^*$.
Note that if $\slX$ is an inverse category, then $\act g{(e^{-1})}=(\acth ge)^{-1}$ for every $g\in G$ and every arrow $e$.
We say that $G$ acts {\em transitively} on $\slX$ if, for any objects $i,j$, there exists $g\in G$ with $j=\acth gi$, and that
$G$ acts on $\slX$ {\em without fixed points} if, for any $g\in G$ and any object $i$, we have $\acth gi=i$ only if $g=1$.
Note that if $G$ acts transitively on $\slX$, then the local monoids $\slX(i,i)\ (i\in\Ob\slX)$ are all isomorphic.

Let $G$ be a group acting on a category $\slX$. 
This action determines a category $\slX/G$ in a natural way:
the objects of $\slX/G$ are the orbits of the objects of $\slX$, and, for every pair $\act Gi,\act Gj$ of objects, 
the $(\act Gi,\act Gj)$-arrows are the orbits of the $(i',j')$-arrows of $\slX$ where $i'\in \act Gi$ and $j'\in \act Gj$.
The product of consecutive arrows $\te,\tf$ is also defined in a natural way, namely, by considering the orbit of a product $ef$
where $e,f$ are consecutive arrows in $\slX$ such that $e\in\te$ and $f\in\tf$.
Note that if $G$ acts transitively on $\slX$, then $\slX/G$ is a one-object category, that is, a monoid.
The properties below are proved in \cite[Propositions 3.11, 3.14]{mar-pin}.

\begin{Res}\label{cat-per-gr}
Let $G$ be a group acting transitively and without fixed points on an inverse category $\slX$.
\begin{enumerate}
\item\label{cat-per-gr1}
The monoid $\slX/G$ is inverse, and it is isomorphic, for every object $i$, to the monoid $(\slX/G)_i$ defined on the set
$\{(e,g): g\in G\ \hbox{and}\ e\in \slX(i,\acth gi)\}$
by the multiplication
\begin{equation*}
(e,g)(f,h)=(e\cdot \act gf, gh).
\end{equation*}
\item\label{cat-per-gr2}
If $\slX$ is connected and it is locally a semilattice, then $\slX/G$ is an $E$-unitary inverse monoid.
Moreover, the greatest group homomorphic image of $\slX/G$  is $G$, and its semilattice of idempotents is isomorphic to
$\slX(i,i)$ for any object $i$.
\item\label{cat-per-gr3}
If $\slX$ is connected, and it is locally in a group variety $\slU$, then $\slX/G$ is a group which is an extension of $\slX(i,i)\in\slU$ by $G$
for any object $i$.
\end{enumerate}
\end{Res}

For our later convenience, note that the inverse of an element can be obtained in $(\slX/G)_i$ in the following manner:
\begin{equation*}
(e,g)^{-1}=(\act{g^{-1}}{e^{-1}},g^{-1}).
\end{equation*}

Notice that if a group $G$ acts on an inverse category transitively and without fixed points, then $\Ob\slX$ is in one-to-one correspondence
with $G$.
In the sequel we consider several categories of this kind which have just $G$ as its set of objects.
For these categories, we identify $\slX/G$ with $(\slX/G)_1$.

To see that each $E$-unitary inverse monoid can be obtained in the way described in Result \ref{cat-per-gr}(\ref{cat-per-gr2}),
let $M$ be an arbitrary $E$-unitary inverse monoid, and denote the group $M/\sigma$ by $G$.
Consider the derived category (see \cite{Tilson}) of the natural homomorphism $\sigma^\natural\colon M\to G$, and denote it $\slI_M$:
its set of objects is $G$, its set of $(i,j)$-arrows is
$$\slI_M(i,j)=\{(i,m,j)\in G\times M\times G: i\cdot m\sigma=j\}\ (i,j\in G),$$
and the product of consecutive arrows $(i,m,j)\in \slI_M(i,j)$ and $(j,n,k)\in \slI_M(j,k)$ is defined by the rule
$$(i,m,j)(j,n,k)=(i,mn,k).$$
It is easy to see that an arrow $(i,m,j)$ is idempotent if and only if $m$ is idempotent, and if this is the case, then $i=j$.
Moreover, we have $(i,m,j)^{-1}=(j,m^{-1},i)$ for every arrow $(i,m,j)$. 
The natural partial order on $\slI_M$ is the following: for any arrows $(i,m,j),(k,n,l)$, we have
$(i,m,j)\le (k,n,l)$ if and only if $i=k,\ j=l$ and $m\le n$.

The group $G$ acts naturally on $\slI_M$ as follows: $\acth gi=gi$ and $\act g{(i,m,j)}=(gi,m,gj)$ for every $g\in G$ and $(i,m,j)\in \Arr{\slI_M}$.

The category $\slI_M$ and the action of $G$ on it has the following properties \cite[Proposition 3.12]{mar-pin}.

\begin{Res}\label{locsl}
The category $\slI_M$ is a connected inverse category which is locally a semilattice.
The group $G$ acts transitively and without fixed points on $\slI_M$, and $M$ is isomorphic to $\slI_M/G$.
\end{Res}

Let $G$ be an $A$-generated group where $A\subseteq G\setminus \{1\}$. 
The {\em Margolis--Meakin expansion $M(G)$ of $G$} is defined in the following way:
consider the set of all pairs $(\rmX,g)$ where $g\in G$ and $\rm X$ is a finite connected subgraph of the Cayley graph $\Gamma$ of $G$
containing the vertices $1$ and $g$, and  define a multiplication on this set by the rule
$$(\rmX,g)(\rmY,h)=(\rmX\cup \acth g{\rmY},gh).$$
Then $M(G)$ is an $A$-generated $E$-unitary inverse monoid with
$\epsilon_{M(G)}\colon A\to M(G),\ a\mapsto (\sp{e_a},a)=(e_a,a)$
(i.e., for brevity, we identify $\sp{e}$ with $e$ for every edge $e$ in $\Gamma$),
where the identity element is $(\emptyset_1,1)$ and $(\rmX,g)^{-1}=(\act{g^{-1}}{\rmX},g^{-1})$ for every $(\rmX,g)\in M(G)$.
In particular, if $G$ is the free $A$-generated group, then $M(G)$ is the free inverse monoid.

By definition, the arrows in $\slI_{M(G)}(i,j)$ are $(i,(\rmX,g),j)$ where $(\rmX,g)\in M(G)$ and $ig=j$ in $G$.
Therefore $\slI_{M(G)}/G=(\slI_{M(G)}/G)_1$ consists of the pairs $((1,(\rmX,g),g),g)$ which can be identified with
$(\rmX,g)$, and this identification is the isomorphism involved in Result \ref{locsl}.
Moreover, notice that the assignment $(i,(\rmX,g),j)\mapsto (i,\acth i{\rmX},j)$ is a bijection from $\slI_{M(G)}(i,j)$
onto the set of all triples $(i,\rmX,j)$ where $\rmX$ is a finite connected subgraph of $\Gamma$ and $i,j\in V_{\rmX}$.
Thus $\slI_{M(G)}$ can be identified with the category where the hom-sets are the latter sets, and the multiplication is the 
following:
$$(i,\rmX,j)(j,\rmY,k)=(i,\rmX\cup\rmY,k).$$

Now we construct the `most general' $A$-generated group which is an extension of a member of a group variety
$\slU$ by a given $A$-generated group $G$, in the form $\slX/G$ where $\slX$ is a category (see \cite{ASzM}).
Consider a group variety $\slU$ and an $A$-generated group $G$. 
Denote the Cayley graph of $G$ by $\Gamma$, and the 
relatively free group in $\slU$
on $E_\Gamma$ by $F_\slU(E_\Gamma)$.
Note that the action of $G$ on $\Gamma$ extends naturally to an action of $G$ on $F_\slU(E_\Gamma)$,
and this defines a semidirect product $F_\slU(E_\Gamma) \rtimes G$.
Any path in $\ol\Gamma$, regarded as a word in $\ol{E_\Gamma}^\ast$, determines an element of $F_\slU(E_\Gamma)$
which is denoted by $[p]_\slU$, as is introduced above.

By \cite{Tilson}, the {\em free $g\slU$-category on $\Gamma$}, denoted by $\cat\slU$, is given as follows: 
its set of objects is $V_\Gamma$, and, for any pair of objects $i,j$, the set of $(i,j)$-arrows is
$$\cat\slU(i,j)=\{(i,[p]_\slU,j): p\hbox{ is a }(i,j)\hbox{-path in }\ol\Gamma\},$$
and the product of consecutive arrows is defined by
$$(i,[p]_\slU,j)(j,[q]_\slU,k)=(i,[pq]_\slU,k).$$
Obviously, the category $\cat\slU$ is a groupoid, and the inverse of an arrow is obtained as follows:
$$(i,[p]_\slU,j)^{-1}=(j,[p]^{-1}_\slU,i)=(j,[p']_\slU,i).$$
Moreover, if $\slU$ is non-trivial, then the map $\epsilon_{\cat\slU}\colon\Gamma\to\cat\slU$,
defined by $e\mapsto (\iota e,[e]_\slU,\tau e)=(\iota e,e,\tau e)$ 
(i.e., as usual, we identify $[e]_\slU$ with $e$ in the free group $F_\slU(E_\Gamma)$)
for every edge $e$ in $\Gamma$, embeds the graph $\Gamma$ into $\cat\slU$, and $\Gamma\epsilon_{\cat\slU}$ generates $\cat\slU$. 

Notice that the action of $G$ on $\Gamma$ extends to an action of $G$ of $\cat\slU$, and this action is transitive and has no fixed points.
Furthermore, $\cat\slU$ is connected since $\Gamma$ is connected.
Thus Result \ref{cat-per-gr}(\ref{cat-per-gr3}) implies that $\cat\slU/G$ is a group which is an extension of a member of $\slU$ by $G$.
What is more, it is straightforward to see by Result \ref{cat-per-gr}(\ref{cat-per-gr1}) that 
the elements of $\cat\slU/G=(\cat\slU/G)_1$ are exactly the pairs $([p]_\slU,g) \in F_\slU(E_\Gamma) \rtimes G$, where $p$ is a $(1,g)$-path in $\ol \Gamma$.
Moreover, $\cat\slU/G$ is generated by the subset $\{(e_a,a\epsilon_G) : a \in A\}$, 
and so it is $A$-generated with $\epsilon_{\cat\slU/G} \colon A \to \cat\slU/G, a \mapsto (e_a, a\epsilon_G)$.
On the other hand, we see that $\cat\slU/G$ is a subgroup in the semidirect product $F_\slU(E_\Gamma) \rtimes G$.
It is well known (cf.\ the Kaloujnine--Krasner theorem) that $\cat\slU/G$ is the `most general' $A$-generated group 
which is an extension of a member of $\slU$ by $G$, that is,
it has the universal property that, for each such extension $K$ with $\epsilon_K\colon A\to K$, there exists a surjective homomorphism 
$\varphi\colon \cat\slU/G\to K$ such that $\epsilon_{\cat\slU/G}\varphi=\epsilon_K$.
For brevity, we denote the group $\cat\slU/G$ later on by $G^\slU$, see \cite{ASzM}.

\subsection*{Dual premorphisms}
For any inverse categories $\slX$ and $\slY$, a graph morphism $\psi \colon \slX \to \slY$ is called a {\em dual premorphism} if 
$1_i\psi=1_{i\psi}$, $(e^{-1})\psi=(e\psi)^{-1}$ and $(ef)\psi\geq e\psi \cdot f\psi$ for any object $i$ and any consecutive arrows $e,f$ in $\slX$. 
In particular, this defines the notion of a dual premorphism between one-object inverse categories, that is, between inverse monoids
(such maps are called dual prehomomorphisms in \cite{lo} and prehomomorphisms in \cite{pet}).

An important class of dual premorphisms from groups to an inverse monoid $M$ is closely related to $F$-inverse covers of $M$, 
as stated in the following well-known result (\cite[Theorem VII.6.11]{pet}):

\begin{Res}\label{dual}
Let $H$ be a group and $M$ be an inverse monoid. If $\psi \colon H \to M$ is a dual premorphism
such that 
\begin{equation}\label{tul}
\hbox{for every }m \in M, \hbox{ there exists } h \in H \hbox{ with } m\leq h\psi, 
\end{equation}
then
$$F=\{(m,h) \in M \times H : m \leq h\psi\}$$
is an inverse submonoid in the direct product $M \times H$, and it is an $F$-inverse cover of $M$ over $H$. 
Conversely, up to isomorphism, every $F$-inverse cover of $M$ over $H$ can be so constructed.
\end{Res}

In the proof of the converse part of Result \ref{dual}, the following dual premorphism $\psi \colon F/\sigma \to M$  is constructed for an inverse monoid $M$, 
an $F$-inverse monoid $F$, and a surjective idempotent-separating homomorphism $\varphi \colon F \to M$:
for every $h \in F/\sigma$, let $h\psi=m_h\varphi$, where $m_h$ denotes the maximum element of the $\sigma$-class $h$. 
It is important to notice that, more generally, this construction gives a dual premorphism with property (\ref{tul})
for any surjective homomorphism $\varphi \colon F \to M$. 
In the sequel, we call this map $\psi$ the \textit{dual premorphism induced by} $\varphi$.

Notice that, for every group $H$ and inverse monoids $M,N$, the product of a dual premorphism $\psi\colon H\to M$ with property (\ref{tul})
and a surjective homomorphism $\varphi\colon M\to N$ is a dual premorphism from $H$ to $N$ with property (\ref{tul}).
As a consequence, notice that if an inverse monoid $M$ has an $F$-inverse cover over a group $H$, then so do its homomorphic images.

\section{Conditions on the existence of $F$-inverse covers}

In this section, the technique introduced in \cite{ASzM} is generalized for a class of $E$-unitary inverse monoids containing all finite ones
and all Margolis--Meakin expansions of $A$-generated groups, 
and necessary and sufficient conditions are provided for any member of this class to have an $F$-inverse cover over a given variety of groups.

First, we define the class of $E$-unitary inverse monoids we intend to consider.
The idea comes from the observation that \cite[Lemma 2.3]{ASzM} remains valid under an assumption weaker than $M$ being $A$-generated
where the elements of $A$ are maximal in $M$ with respect to the natural partial order.
We introduce the appropriate notion more generally for inverse categories.

Let $\slX$ be an inverse category and $\Delta$ an arbitrary graph. 
We say that $\slX$ is {\em quasi-$\Delta$-generated} if a graph morphism $\epsilon_\slX\colon \Delta\to \slX$ is given such that
the subgraph $\Delta\epsilon_\slX\cup \rmE(\slX)$ generates $\slX$, where $\rmE(\slX)$ is the subgraph of the idempotents of $\slX$. 
Clearly, a $\Delta$-generated inverse category is quasi-$\Delta$-generated.
Furthermore, notice that a groupoid is quasi-$\Delta$-generated if and only if it is $\Delta$-generated. 
If $\epsilon_\slX$ is injective, then we might assume that $\Delta$ is a subgraph in $\slX$, i.e., $\epsilon_\slX$ is the inclusion graph morphism 
$\Delta\to \slX$.

A dual premorphism $\psi\colon \slY\to \slX$ between quasi-$\Delta$-generated inverse categories is called {\em canonical} if $\epsilon_\slY\psi=\epsilon_\slX$.
Note that, in this case, if $\epsilon_\slX$ is an inclusion, then $\epsilon_\slY$ is necessarily injective, and so it also can be chosen to be an inclusion.
However, if $\epsilon_\slY$ is injective (in particular, an inclusion), then $\epsilon_\slX$ need not be injective, and so one cannot suppose in general that
$\epsilon_\slX$ is an inclusion.

In particular, if $\slX,\slY$ are one-object inverse categories, that is, inverse monoids, and $\Delta$ is a one-vertex graph, that is, a set, then this defines 
a quasi-$A$-generated inverse monoid and a canonical dual premorphism between inverse monoids. 
We also see that a group is quasi-$A$-generated if and only if $A$-generated. 

An inverse monoid $M$ is called {\em finite-above} if the set $m^\omega=\{n\in M:n\ge m\}$ is finite for every $m\in M$.
For example, finite inverse monoids and the Margolis--Meakin expansions of $A$-generated groups are finite-above.
The class we investigate in this section is that of all finite-above $E$-unitary inverse monoids.

Notice that if $M$ is a finite-above inverse monoid, then, for every element $m\in M$, there exists $m'\in M$ 
such that $m'\ge m$ and $m'$ is maximal in $M$ with respect to the natural partial order. 
Denoting by $\maxM$ the set of all elements of $M$ distinct from $1$ which are maximal with respect to the natural partial order, 
we obtain that $M$ is quasi-$\maxM$-generated.
Hence the following is straightforward.

\begin{Lem}\label{fa-qAg}
Every finite-above inverse monoid is quasi-$A$-generated for some $A\subseteq\maxM$.
\end{Lem}

What is more, the following lemma shows that each quasi-generating set of a finite-above inverse monoid 
can be replaced in a natural way by one contained in $\maxM$.
As usual, the set of idempotents $E(M)$ of $M$ is simply denoted by $E$.
Note that if $A\subseteq\maxM$, then $A\cap E=\emptyset$.
Here and later on, we need the following notation.
If $M$ is quasi-$A$-generated and $w$ is a word in $\ol{A\cup E}^*$, then the word in $\ol{A\setminus E}^*\subseteq \ol A^*$
obtained from $w$ by deleting all letters from $\ol E$ is denoted by $w^-$.
Obviously, we have $[w]_M\le [w^-]_M$.

\begin{Lem}\label{mx}
Let $M$ be a finite-above inverse monoid, and assume that $A\subseteq M$ is a quasi-generating set in $M$.
For every $a\in A$, let us choose and fix a maximal element $\ta$ such that $a\le \ta$.
Then $\tA=\{\ta:a\in A\}\setminus\{1\}$ is a quasi-generating set in $M$ such that $\tA\subseteq\maxM$.
\end{Lem}

\begin{Proof}
Since $A$ is a quasi-generating set, for every $m\in M$, there exists a word $w\in \ol{A\cup E}^*$ 
such that $m=[w]_M$, whence $m\le [w^-]_M$ follows.
Moreover, the word $\tu$ obtained from $u=w^-$ by substituting $\ta$ for every $a\in A\setminus E$ has the property that
$[u]_M\le [\tu]_M$, and so $m\le [\tu]_M$ holds. 
Thus $m$ belongs to the inverse submonoid of $M$
generated by $\tA\cup E$.
\end{Proof}

This observation establishes that, within the class of finite-above inverse monoids, it is natural to restrict our consideration to quasi-generating sets 
contained in $\maxM$.
Now we present a statement on the $E$-unitary covers of finite-above inverse monoids.

\begin{Lem}
Let $M$ be an inverse monoid.
\begin{enumerate}
\item If $M$ is finite-above, then so are its $E$-unitary covers.\label{E-fa}
\item If $M$ is quasi-$A$-generated for some $A\subseteq \maxM$, then every $E$-unitary cover of $M$ contains a quasi-$A$-generated inverse submonoid $T$
with $A\epsilon_T\subseteq\max T^-$ such that $T$ is an $E$-unitary cover of $M$.\label{E-mg}
\end{enumerate}
\end{Lem}

\begin{Proof}
Let $U$ be any $E$-unitary cover of $M$, and let $\varphi\colon U\to M$ be an idempotent separating and surjective homomorphism.

(\ref{E-fa})
Since $\varphi$ is order preserving, we have $t^\omega\varphi\subseteq (t\varphi)^\omega$ for every $t\in U$, and
the latter set is finite by assumption.
To complete the proof, we verify that 
$\varphi|_{t^\omega}\ (t\in U)$ is injective.
Let $t\in U$ and $y,y_1\in t^\omega$ such that $y\varphi=y_1\varphi$.
This equality implies $yy^{-1}=y_1y_1^{-1}$,  since $\varphi$ is idempotent separating.
Moreover, the relation $y,y_1\ge t$ implies $y\,\sigma\,t\,\sigma\,y_1$, and so
we deduce $y=y_1$, since $U$ is $E$-unitary.

(\ref{E-mg})
For every $a\in A$, let us choose and fix an element $u_a\in U$ such that $u_a\varphi=a$, consider the inverse submonoid
$T$ of $U$ generated by the set $\{u_a:a\in A\}\cup E(U)$, and put $\epsilon_T\colon A\to T,\ a\mapsto u_a$ which is clearly injective.
Obviously, $T$ is a quasi-$A$-generated $E$-unitary inverse monoid, and the restriction $\varphi|_T\colon T\to M$ of $\varphi$ 
is an idempotent separating and surjective homomorphism.
It remains to verify that $A\epsilon_T\subseteq\max T^-$.
Observe that an element $m\in M$ is maximal if and only if the set $m^\omega$ is a singleton, and similarly for $T$.
Thus the last part of the proof of (\ref{E-fa}) shows that $A\epsilon_T\subseteq\max T$.
Since, for every $a\in A$, the relation $a\not=1$ implies $u_a\not= 1$, the proof is complete.
\end{Proof}

This implies the following statement.

\begin{Cor}
Each quasi-$A$-generated finite-above inverse monoid $M$ with $A\subseteq\maxM$ has an $E$-unitary cover with the same properties.
\end{Cor}

This shows that the study of the $F$-inverse covers of finite-above inverse monoids can be reduced to the study of the $F$-inverse covers of finite-above $E$-unitary inverse monoids
in the same way as in the case of finite inverse monoids generated by their maximal elements, see \cite{ASzM}.
Furthermore, the fundamental observations \cite[Lemmas 2.3 and 2.4]{ASzM} can be easily adapted to quasi-$A$-generated finite-above inverse monoids.

\begin{Lem}
Let $H$ be an $A$-generated group and $M$ a quasi-$A$-generated inverse monoid. Then any canonical dual premorphism from $H$ to $M$ has property {\rm (\ref{tul})}.
\end{Lem}

\begin{Proof}
Consider a canonical dual premorphism $\psi \colon H \to M$, and let $m\in M$.
Since $M$ is quasi-$A$-generated, we have $m=[w]_M$ for some $w\in \ol{A\cup E}^*$, and so
$m\le [w^-]_M$ where $w^-\in {\ol A}^*$.
Since $\psi$ is a canonical dual premorphism, we obtain that $[w^-]_H\psi \geq [w^-]_M \geq m$.
\end{Proof}

\begin{Lem}\label{subgroup}
Let $M$ be a quasi-$A$-generated inverse monoid such that $A \subseteq \maxM$.
If $M$ has an $F$-inverse cover over a group $H$, then 
there exists an $A$-generated subgroup $H'$ of $H$ and a canonical dual premorphism from $H'$ to $M$. 
\end{Lem}

\begin{Proof}
Let $F$ be an $F$-inverse monoid and $\varphi\colon F\to M$ a surjective homomorphism. 
Put $H=F/\sigma$, and consider the dual premorphism $\psi\colon H\to M,\ h\mapsto m_h\varphi$ induced by $\varphi$.
Since $\psi$ has property (\ref{tul}), 
for any $a \in A$, there exists $h_a \in H$ such that $a \leq h_a\psi$. 
However, since $a$ is maximal in $M$, this implies $a =h_a\psi$. 
Now let $H'$ be the subgroup of $H$ generated by $\{h_a : a\in A\}$. 
Then the restriction $\psi|_{H'}\colon H'\to M$ of $\psi$ is obviously a dual premorphism. 
Moreover, the subgroup $H'$ is $A$-generated with $\epsilon_{H'} \colon A \to H', a \mapsto h_a$, so $\psi|_{H'}$ is also canonical.
\end{Proof}

So far, the question of whether a finite-above inverse monoid $M$ has an $F$-inverse cover over the class of groups $\slC$ closed under taking
subgroups has been reduced to the question of whether there is a canonical dual premorphism from an $A$-generated group in $\slC$ to $M$,
where $A\subseteq \maxM$ is a quasi-generating set in $M$. 
The answer to this question does not depend on the choice of $A$.

Let $M$ be a quasi-$A$-generated inverse monoid with $A\subseteq \maxM$, $H$ an $A$-generated group in $\slC$, and 
let $\psi\colon H\to M$ be a canonical dual premorphism.
Denote the $A$-generated group $M/\sigma$ by $G$,
and note that $\sigma^\natural \colon M \to G$ is canonical. 
The product  $\kappa=\psi\sigma^\natural$ is a canonical dual premorphism from $H$ to $G$.
However, a dual premorphism between groups is necessarily a homomorphism. 
Consequently, $\kappa \colon H \to G$ is a canonical, and therefore surjective, homomorphism. 
Hence $H$ is an $A$-generated extension of a group $N$ by the $A$-generated group $G$.
If $F$ is an $F$-inverse cover of $M$ over $H$ then, 
to simplify our terminology, we also say that $F$ is an $F$-inverse cover of $M$ {\em via} $N$ or {\em via} a class $\slD$ of groups if $N\in\slD$.
If we are only interested in whether $M$ has an $F$-inverse cover via a member of a given group variety $\slU$, 
then we may replace $H$ by the `most general' $A$-generated extension $G^\slU$ of a member of $\slU$ by $G$.
Thus Lemma \ref{subgroup} implies the following assertion.

\begin{Prop}\label{gen2can}
Let $M$ be a quasi-$A$-generated inverse monoid
with $A\subseteq\maxM$, 
put $G=M/\sigma$, and let $\slU$ be a group variety. 
Then $M$ has an $F$-inverse cover via the group variety $\slU$ if and only if there exists a canonical dual premorphism $G^\slU \to M$.
\end{Prop}

Therefore our question to be studied is reduced to the question of whether there exists a canonical dual premorphism $G^\slU \to M$ with $G=M/\sigma$ 
for a given group variety $\slU$ and for a given quasi-$A$-generated inverse monoid $M$ with $A\subseteq \maxM$.
In the sequel, we deal with this question in the case where $M$ is finite-above and $E$-unitary.

Let $M$ be an $E$-unitary inverse monoid, denote $M/\sigma$ by $G$, and consider the inverse category $\slI_M$ acted upon by $G$.
Given a path $p=e_1e_2\cdots e_n$ in $\ol{\slI_M}$ where $e_j=(\iota e_j,m_j,\tau e_j)$ for every $j\ (j=1,2,\ldots,n)$, consider the word $w=m_1m_2\cdots m_n \in \ol M^\ast$
determined by the labels of the arrows in $p$, and
let us assign an element of $M$ to the path $p$ by defining $\lbl(p)=[w]_M$.
Notice that, for every path $p$, we have $\lbl(p)=\lbl(pp'p)$, and 
$\lbl(p)$ is just the label of the arrow $p\varphi$, where $\varphi\colon \ol{\slI_M}^*\to \slI_M$ is the unique category
morphism such that $e\varphi=e$ and $e'\varphi=e^{-1}$ for every $e\in\Arr\slI_M$. 
Since the local monoids of the category $\slI_M$  are semilattices by Result \ref{locsl},
the following property follows from \cite[Lemma 2.6]{mar-mea} (see also \cite[Chapter VII]{eil} and \cite[Section 12]{Tilson}).

\begin{Lem}
\label{spanning}
For any coterminal paths $p,q$ in $\ol{\slI_M}$, if $\sp{p}=\sp{q}$, then $\lbl(p)=\lbl(q)$.
\end{Lem}

This allows us to assign an element of $M$
to any birooted finite connected subgraph:
if $\rmX$ is a finite connected subgraph in $\slI_M$ and $i,j\in V_\rmX$, then
let $\lbl_{(i,j)}(\rmX)$ be $\lbl(p)$, where $p$ is an $(i,j)$-path in $\ol{\slI_M}$ with $\sp{p}=\rmX$.

Now assume that $M$ is a quasi-$A$-generated $E$-unitary inverse monoid with $A\subseteq \maxM$, and recall that in this case, $G=M/\sigma$ is an $A$-generated group.
Based on the ideas in \cite{mar-mea}, we now give a model for $\slI_M$  as a quasi-$\Gamma$-generated inverse category
where $\Gamma$ is the Cayley graph of $G$.
Choose and fix a subset $I$ of $E$ such that $A\cup I$ generates $M$. 
In particular, if $M$ is $A$-generated, then $I$ can be chosen to be empty.
Consider the subgraphs $\Gamma$ and $\Gamma^I$ of $\slI_M$ consisting of all edges with labels from $A$ and from $A\cup I$, respectively.
Notice that $\Gamma$ is, in fact, the Cayley graph of the $A$-generated group $G$, and $\Gamma^I$ is obtained from $\Gamma$ by adding loops to it (with labels from $I$).

We are going to introduce a closure operator on the set $\Sub(\Gamma^I)$ of all subgraphs of $\Gamma^I$.
We need to make a few observations before.

\begin{Lem}
\label{largestlabel}
Let $\rmX,\rmY$ be finite connected subgraphs in $\Gamma^I$, and let $i,j \in V_\rmX \cap V_\rmY$.
If $\lbl_{(i,j)}(\rmX)\le \lbl_{(i,j)}(\rmY)$, then $\lbl_{(i,j)}(\rmX)=\lbl_{(i,j)}(\rmX\cup \rmY)$.
\end{Lem}

\begin{Proof}
Let $r$ and $s$ be arbitrary $(i,j)$-paths spanning $\rmX$ and $\rmY$, respectively. 
Then $rr's$ is an $(i,j)$-path spanning $\rmX \cup \rmY$. 
According to the assumption, $\lbl(r)\le \lbl(s)$, so $\lbl(rr's)=\lbl(r)$. 
\end{Proof}

\begin{Lem}
\label{pointlesslabel}
Let $\rmX,\rmY$ be finite connected subgraphs in $\Gamma^I$, and let $i,j \in V_\rmX \cap V_\rmY$. 
If $\lbl_{(i,j)}(\rmX)\le\lbl_{(i,j)}(\rmY)$, then $\lbl_{(k,l)}(\rmX)\le\lbl_{(k,l)}(\rmY)$ for every $k,l \in V_\rmX \cap V_\rmY$.
\end{Lem}

\begin{Proof}
Let $r$ and $s$ be $(i,j)$-paths spanning $\rmX$ and $\rmY$, respectively, and let $p_1$ and $q_1$ be $(k,i)$-paths in $\rmX$ and $\rmY$, and 
let $p_2$ and $q_2$ be $(j,l)$-paths in $\rmX$ and $\rmY$, respectively. 
Then $p_1r p_2$ and $q_1sq_2$ are $(k,l)$-paths spanning $\rmX$ and $\rmY$, respectively. 
Therefore, by applying Lemmas \ref{spanning} and \ref{largestlabel}, we obtain that
\begin{eqnarray*}
\lbl_{(k,l)}(\rmX)&\!\!\!=\!\!\!&\lbl(p_1rp_2)=\lbl(p_1)\lbl_{(i,j)}(\rmX)\lbl(p_2)=\lbl(p_1)\lbl_{(i,j)}(\rmX\cup \rmY)\lbl(p_2)\cr
                  &\!\!\!=\!\!\!&\lbl(p_1)\lbl(rr's)\lbl(p_2)=\lbl(p_1rr'sp_2)=\lbl(q_1rr'sq_2)\cr
                  &\!\!\!\le\!\!\!& \lbl(q_1sq_2)=\lbl_{(k,l)}(\rmY).
\end{eqnarray*}
\end{Proof}

Given a finite connected subgraph $\rmX$ in $\Gamma^I$ with vertices $i,j\in V_\rmX$, consider the subgraph 
\begin{align*}
\cl \rmX=\bigcup\{\rmY\in\Sub(\Gamma^I):\ &\rmY\hbox{ is finite and connected, } i,j\in V_\rmY,\\
& \hbox{and }\lbl_{(i,j)}(\rmY)\ge \lbl_{(i,j)}(\rmX)\}
\end{align*}
of $\Gamma^I$ which is clearly connected. 
Note that, by Lemma \ref{pointlesslabel}, the graph $\cl \rmX$ is independent of the choice of $i,j$.
Moreover, by Lemma \ref{largestlabel},  
the same subgraph is obtained if the relation `$\ge$' is replaced by `$=$' in the definition of $\cl \rmX$.
More generally, for any $\rmX\in\Sub(\Gamma^I)$, let us define the subgraph $\cl \rmX$ in the following manner:
$$\cl \rmX=\bigcup\{\cl \rmY:\rmY\hbox{ is a finite and connected subgraph of } \rmX\}.$$

It is routine to check that $\rmX\to \cl \rmX$ is a closure operator on $\Sub(\Gamma^I)$, that is,
$\rmX\subseteq \cl \rmX$, $\cl{(\cl \rmX)}=\cl \rmX$, and $\rmX\subseteq\rmX_1$ implies $\cl \rmX\subseteq \cl{\rmX_1}$
for any $\rmX,\rmX_1\in\Sub(\Gamma^I)$.
As usual, a subgraph $\rmX$ of $\Gamma^I$ is said to be {\em closed} if $\rmX=\cl \rmX$.
Note that, in particular, we have
$$\cl{\emptyset}_i=\bigcup\left\{\sp{h}:h\ \hbox{is an}\ i\hbox{-cycle in } \Gamma^I \hbox{ such that}\ \lbl(h)=1\right\},$$ 
and so $\emptyset_i$ is closed if and only if there is 
no $a\in A$ such that $a\,\mathcal{R}\, 1$ or $a\,\mathcal{L}\, 1$.
Furthermore, we have $\cl{\rmX}\supseteq \cl{\emptyset}_i$ for every $\rmX\in\Sub(\Gamma^I)$ and $i\in V_{\cl\rmX}$. 
In particular, we see that the closure of a finite subgraph need not be finite.
For example, if $M$ is the bicyclic inverse monoid generated by 
$A=\{a\}$ where $aa^{-1}=1$, then $a$ is a maximal element in $M$,
$M/\sigma$ is the infinite cyclic group generated by $a\sigma$, and we have 
$\cl{\emptyset}_1=\{((a\sigma)^n,a,(a\sigma)^{n+1}): n\in\mathbb{N}_0\}$.

Denote the set of all closed subgraphs of $\Gamma^I$ by $\ClSub(\Gamma^I)$, and its subset consisting of the closures of all finite connected subgraphs by
$\ClSubfc(\Gamma^I)$. 
Moreover, for any family $\rmX_j\ (j\in J)$ of subgraphs of $\Gamma^I$, define $\bigvee_{j\in J}\rmX_j=\cl{\bigl(\bigcup_{j\in J}\rmX_j\bigr)}$.
The following lemmas formulate important properties of closed subgraphs which can be easily checked.

\begin{Lem}\label{cl}
For every quasi-$A$-generated $E$-unitary inverse monoid $M$ with $A\subseteq \maxM$, the following statements hold.
\begin{enumerate}
\item Each component of a closed subgraph is closed.
\item The partially ordered set $(\ClSub(\Gamma^I);\subseteq)$ forms a complete lattice with respect to the usual intersection and 
the operation $\bigvee$ defined above.
\item For any $\rmX,\rmY\in \ClSubfc(\Gamma^I)$ with $V_\rmX\cap V_\rmY\not=\emptyset$, we have  $\rmX\vee\rmY\in \ClSubfc(\Gamma^I)$.
\item \label{cl4}
For any finite connected subgraph in $\Gamma^I$ and for any $g\in G$, we have $\acth g{(\cl{\rmX})}=\cl{(\acth g{\rmX})}$.
Consequently, the action of $G$ on $\Sub(\Gamma^I)$ restricts to an action on $\ClSub(\Gamma^I)$ and to an action on $\ClSubfc(\Gamma^I)$,
respectively.
\end{enumerate}
\end{Lem}

Now we prove that the descending chain condition holds for $\ClSubfc(\Gamma^I)$ if $M$ is finite-above.

\begin{Lem}\label{fc-acc}
If $M$ is a quasi-$A$-generated finite-above $E$-unitary inverse monoid with $A\subseteq \maxM$, then, for every
$\rmX\in \ClSubfc(\Gamma^I)$ and $i\in V_\rmX$, there are only finitely many closed connected subgraphs in $\rmX$ containing the vertex $i$, 
and all belong to $\ClSubfc(\Gamma^I)$.
\end{Lem}

\begin{Proof}
Let $\rmX\in\ClSubfc(\Gamma^I)$, whence $\rmX=\cl \rmY$ for some finite connected subgraph $\rmY$, and let $i\in V_\rmY$.
If $\rmZ$ is any finite connected subgraph such that $\rmX\supseteq \cl \rmZ$ and $i\in V_\rmZ$, then 
$\lbl_{(i,i)}(\rmY)\le \lbl_{(i,i)}(\rmZ)$.
Since $M$ is finite-above, the set 
$\Lambda=\{\rmX_0\in\ClSubfc(\Gamma^I):\rmX_0\subseteq \rmX \hbox{ and } i\in V_{\rmX_0}\}$
is finite.
If $\rmX_1\in\ClSub(\Gamma^I)$ is connected with $\rmX_1\subseteq \rmX$ and $i\in V_{\rmX_1}$,
then, by definition, $\rmX_1$ is a join of a subset of the finite set $\Lambda$ which is closed under $\vee$.
Hence it follows that $\rmX_1$ belongs to $\Lambda$, i.e., $\rmX_1\in\ClSubfc(\Gamma^I)$.
\end{Proof}

Now we define an inverse category $\xc$ in the following way: its set of objects is $G$, its set of $(i,j)$-arrows $(i,j \in G)$ is
$$\xc(i,j)=\{(i, \rmX, j): \rmX \in \ClSubfc(\Gamma) \hbox{ and } i,j \in V_\rmX\},$$ 
and the product of two consecutive arrows is defined by
$$(i, \rmX, j)(j, \rmY, k)=(i, \rmX \vee \rmY, k).$$
It can be checked directly (see also \cite{mar-mea}) that $\xc\to \slI_M,\ (i,\rmX,j)\mapsto(i,\lbl_{(i,j)}(\rmX),j)$
is a category isomorphism. 
Hence $\xc$ is an inverse category with $(i, \rmX, j)^{-1}=(j, \rmX, i)$, 
it is locally a semilattice, and  
the natural partial order on it is the following: $(i, \rmX, j) \leq (k, \rmY, l)$ if and only if $i=k, j=l$ and $\rmX \supseteq \rmY$.
Moreover, the group $G$ acts on it by the rule $\act g{(i,\rmX,j)}=(gi,\acth g{\rmX},gj)$ transitively and without fixed points.
The inverse category $\xc$ is $\Gamma^I$-generated with $\epsilon^I_{\xc} \colon \Gamma^I \to \xc, e \mapsto (\iota e, \cl e, \tau e)$. 
Therefore $\xc$ is also quasi-$\Gamma$-generated with $\epsilon_{\xc}=\epsilon^I_{\xc}|_\Gamma \colon \Gamma \to \xc$.
By Results \ref{cat-per-gr} and \ref{locsl}, hence we deduce the following proposition.

\begin{Prop}\label{xc2M}
\begin{enumerate}
\item The $E$-unitary inverse monoid $\xc/G$ can be described, up to isomorphism, in the following way: its underlying set is
$$\xc/G=\{(\rmX,g): \rmX \in \ClSubfc(\Gamma^I),\ 1, g\in V_\rmX\},$$
and the multiplication is defined by
$$(\rmX,g)(\rmY,h)=(\rmX \vee \acth g{\rmY}, gh).$$
\item The monoid $\xc/G$ is quasi-$A$-generated with 
$$\epsilon_{\xc/G}\colon A\to \xc/G,\quad a\mapsto (\cl{e}_a,a\sigma).$$
\item \label{xc2M3}
The map $\varphi \colon \xc/G \to M,\ (\rmX, g) \mapsto \lbl_{(1,g)}(\rmX)$ is a canonical isomorphism.
\end{enumerate}
\end{Prop}

\begin{Rem}
Proposition \ref{xc2M} provides a representation of $M$ as a $P$-semigroup. 
The McAlister triple involved consists of $G$, the partially ordered set $(\ClSubfc(\Gamma^I);\subseteq)$ and its
order ideal and subsemilattice $(\{\rmX \in \ClSubfc(\Gamma^I),\ 1\in V_\rmX\};\vee)$.
\end{Rem}

Notice that if we apply the construction before Proposition \ref{xc2M} for $M$ being the Margolis--Meakin expansion $M(G)$ 
of an $A$-generated group $G$ with $A\subseteq G\setminus\{1\}$, then $\Gamma^I=\Gamma$, the Cayley graph of $G$, 
the closure operator $\rmX\to \cl \rmX$ is identical on $\Sub(\Gamma)$, and the operation $\vee$ coincides with the usual $\cup$.
Thus the category $\xc$ is just the category isomorphic to $\slI_{M(G)}$ which is presented after Result \ref{locsl},
and the map $\varphi$ given in the last statement of the proposition is, in fact, identical.

The goal of this section is to give equivalent conditions for the existence of a canonical dual premorphism $G^\slU \to M$. 
The previous proposition reformulates it by replacing $M$ with $\xc/G$.
Since $G^\slU=\cat\slU/G$, it is natural to study the connection between the canonical dual premorphisms $\cat\slU/G\to \xc/G$
and the canonical dual premorphisms $\cat\slU\to \xc$.
As one expects, there is a natural correspondence between these formulated in
the next lemma in a more general setting.
The proof is straightforward, it is left to the reader.

\begin{Lem}\label{two-dual}
Let $\Delta$ be any graph, and let $\slY$ be a $\Delta$-generated, and $\slX$ a quasi-$\Delta$-generated inverse category containing $\Delta$.
Suppose that $G$ is a group acting on both $\slX$ and $\slY$ transitively and without fixed points in a way that $\Delta$ is invariant with respect 
to both actions, and the two actions coincide on $\Delta$.
Let $i$ be a vertex in $\Delta$.
\begin{enumerate}
\item \label{two-dual1}
We have $\Ob\slX=V_\Delta=\Ob\slY$, and so the actions of $G$ on $\Ob\slX$ and $\Ob\slY$ coincide.
\item\label{two-dual2}
The inverse monoid $\slY_i$ is 
$\Delta(i,-)$-generated, and
the inverse monoid $\slX_i$ is 
quasi-$\Delta(i,-)$-generated
with the maps 
$$\epsilon_{\slY_i}\colon \Delta(i,-)\to \slY_i,\ e\mapsto (e,g),\quad \hbox{provided } e\in\slY(i,\acth gi),$$
and
$$\epsilon_{\slX_i}\colon \Delta(i,-)\to \slX_i,\ e\mapsto (e,g),\quad \hbox{provided } e\in\slX(i,\acth gi),$$
respectively.
\item \label{two-dual3}
If $\Psi\colon \slY\to \slX$ is a canonical dual premorphism such that 
\begin{equation}\label{commut}
(\act gy)\Psi=\act g{(y\Psi)}\quad \hbox{for every } g\in G \hbox{ and } y\in\Arr\slY, 
\end{equation}
then $\iota (y\Psi)=\iota y$, $\tau (y\Psi)=\tau y$, and 
the map $\psi\colon \slY_i\to \slX_i,\ (e,g)\mapsto (e\Psi,g)$
is a canonical dual premorphism.
\item \label{two-dual4}
If $\psi\colon \slY_i\to \slX_i$ is a canonical dual premorphism and $(e,g)\psi=(\tilde e,\tilde g)$ for some $(e,g)\in\slY_i$ and $(\tilde e,\tilde g)\in\slX_i$,
then $g=\tilde g$, $\iota e=\iota \tilde e$ and $\tau e=\tau \tilde e$.
Thus a graph morphism 
$\Psi\colon \slY\to \slX$ 
can be defined such that, for any arrow $y\in\slY(\acth gi,\acth hi)$, we set $y\Psi$ to be the unique arrow $x\in\slX(\acth gi,\acth hi)$ such that
$(\act{g^{-1}}y,g^{-1}h)\psi=(\act{g^{-1}}x,g^{-1}h)$.
This $\Psi$ is a canonical dual premorphism satisfying {\rm (\ref{commut})}.
\end{enumerate}
\end{Lem}

From now on, let $M$ be a quasi-$A$-generated finite-above $E$-unitary inverse monoid with $A\subseteq\maxM$, and let $\slU$ be an arbitrary group variety. 
Motivated by Lemma \ref{two-dual}, we intend to find a necessary and sufficient condition in order that a canonical dual premorphism $\cat\slU\to \xc$
exists fulfilling condition (\ref{commut}).

We are going to assign two series of subgraphs of $\Gamma^I$ to any arrow $x$ of $\cat \slU$. 
Let
$$\Cc{0}(x)=\bigcap\{\spc{p}: p \hbox{ is a } (\iota x, \tau x)\hbox{-path in } \ol{\Gamma}
\hbox{ such that }x=(\iota x, [p]_\slU, \tau x)\},$$
and let $\Pc{0}(x)$ be the component of $\Cc{0}(x)$ containing $\iota x$. 
Suppose that, for some $n\ (n\ge 0)$, the subgraphs $\Cc{n}(x)$ and $\Pc{n}(x)$ are defined for every arrow $x$ of $\cat \slU$. Then let
\begin{align*}
\Cc{n+1}(x)=\bigcap\{\Pc{n}(&x_1) \vee \cdots \vee \Pc{n}(x_k): k\in \mathbb{N}_0,\ x_1, \ldots, x_k \in \cat\slU  \\
&\hbox{ are consecutive arrows, and }x=x_1\cdots x_k\},
\end{align*}
and again, let $\Pc{n+1}(x)$ be the component of $\Cc{n+1}(x)$ containing $\iota x$. 
Applying Lemma \ref{cl} we see that, 
for every $n$, the subgraph $\Pc{n}(x)$ of $\Gamma^I$ is a component of an intersection of closed subgraphs, so $\Pc{n}(x) \in \ClSub(\Gamma^I)$ and is connected. 
Also, $\Pc{n}(x)$ contains $\iota x$ for all $n$.
Moreover, observe that
$$\Cc{0}(x) \supseteq \Pc{0}(x) \supseteq \cdots \supseteq \Cc{n}(x) \supseteq \Pc{n}(x) \supseteq \Cc{n+1}(x) \supseteq \Pc{n+1}(x) \supseteq \cdots$$
for all $x$ and $n$. 
By Lemma \ref{fc-acc} we deduce that, for every $x$, all these subgraphs belong to $\ClSubfc(\Gamma^I)$, and there exists $n_x\in\mathbb{N}_0$ such that
$\Pc{n_x}(x)=\Pc{n_x+k}(x)$ for every $k\in\mathbb{N}_0$.
For brevity, denote $\Pc{n_x}(x)$ by $\Pc{}(x)$.
Furthermore, for any consecutive arrows $x$ and $y$, we have
$$\Pc{n+1}(xy) \subseteq \Cc{n+1}(xy) \subseteq \Pc{n}(x) \vee \Pc{n}(y),$$
and so
\begin{equation*}
\Pc{}(xy) \subseteq \Pc{}(x) \vee \Pc{}(y)
\end{equation*}
is implied.

\begin{Prop}\label{nemszakad}
There exists a canonical dual premorphism $\psi \colon \cat\slU \to \xc$ if and only if $\Pc{n}(x)$ contains  
$\tau x$ for every $n\in\mathbb{N}_0$ and for every $x \in \cat\slU$, or, equivalently, if and only if $\Pc{}(x)$ contains 
$\tau x$ for every $x \in \cat\slU$.
\end{Prop}

\begin{Proof}
Let $\psi \colon \cat\slU \to \xc$ be a canonical dual premorphism. 
We denote the middle entry of $x\psi$ by $\mu(x\psi)$, 
which belongs to $\ClSubfc(\Gamma^I)$ and contains $\iota x$ and $\tau x$. 
The fact that $\psi$ is a dual premorphism means that $\mu((xy)\psi)\subseteq \mu(x\psi) \vee \mu(y\psi)$.
Moreover, $\psi$ is canonical, therefore we have $(\iota e, [e]_\slU, \tau e)\psi=(\iota e, \cl{e}, \tau e)$ for every $e \in E_\Gamma$. 
Hence for an arbitrary representation of an arrow $x=(\iota x, [p]_\slU, \tau x)$, where $p=e_1\cdots e_n$ is a 
$(\iota x,\tau x)$-path in $\ol\Gamma$ and 
$e_1,\ldots,e_n\in E_{\ol\Gamma}$, we have 
\begin{eqnarray*}
\mu(x\psi) &\!\!\!\subseteq\!\!\!& \mu((\iota e_1, [e_1]_\slU, \tau e_1)\psi) \vee \cdots \vee \mu((\iota e_n, [e_n]_\slU, \tau e_n)\psi)\\
&\!\!\!=\!\!\!& \cl{e}_1 \vee \cdots \vee \cl{e}_n=\spc{p},
\end{eqnarray*}
which implies $\mu(x\psi) \subseteq \Cc{0}(x)$. 
Since $\mu(x\psi )$ is connected and contains $\iota x$, $\mu( x\psi) \subseteq \Pc{0}(x)$, and this implies $\tau x \in \Pc{0}(x)$.

Now suppose $n \geq 0$ and $\mu (y\psi) \subseteq \Pc{n}(y)$ for any arrow $y$. 
Let $x=x_1\cdots x_k$ be an arbitrary decomposition in $\cat\slU$. 
Then 
$$\mu(x\psi) \subseteq \mu(x_1\psi) \vee \cdots \vee \mu(x_k\psi) \subseteq \Pc{n}(x_1) \vee \cdots \vee \Pc{n}(x_k)$$
holds, whence $\mu(x\psi) \subseteq \Cc{n+1}(x)$. 
As before, $\mu(x\psi)$ is connencted and contains both $\iota x$ and $\tau x$, so we see that $\mu(x\psi) \subseteq \Pc{n+1}(x)$ and $\tau x \in \Pc{n+1}(x)$. 
This proves the `only if' part of the statement.

For the converse, suppose that for any arrow $x$ in $\cat\slU$, we have $\tau x \in \Pc{n}(x)$ for all $n \in \mathbb N_0$. 
We have seen above that $\Pc{}(x)\in\ClSubfc(\Gamma^I)$, and $\Pc{}(xy) \subseteq \Pc{}(x) \vee \Pc{}(y)$ for any arrows $x,y$.
Furthermore, the equality $\Pc{}(x)=\Pc{}(x^{-1})$ can be easily checked for all arrows $x$ by definition. 
Now consider the map $\Pc{}$ which assigns the arrow $(\iota x, \Pc{}(x), \tau x)$ of $\xc$ to the arrow $x$ of $\cat\slU$. 
By the previous observations, this is a dual premorphism from $\cat\slU$ to $\xc$, and the image of $(\iota e, [e]_\slU, \tau e)$ is $(\iota e, \cl{e}, \tau e)$, 
hence it is also canonical.
\end{Proof}

The canonical dual premorphism $\Pc{}$ constructed in the previous proof has property {\rm (\ref{commut})}.

\begin{Lem}\label{aut}
For every $g\in G$ and for any arrow $x$ of $\cat\slU$, we have 
$\Pc{}(\acth g x)=\acth g {\Pc{}(x)}$.
\end{Lem}

\begin{Proof}
One can see by definition that $\Cc{0}(\acth g x)=\acth g {\Cc{0}(x)}$ for all $x \in \cat\slU$, and so $\Pc{0}(\acth g x)=\acth g {\Pc{0}(x)}$ also holds. 
By making use of Lemma \ref{cl}(\ref{cl4}), an easy induction shows that $\Cc{n}(\acth g x)=\acth g {\Cc{n}(x)}$ and $\Pc{n}(\acth g x)=\acth g {\Pc{n}(x)}$ for all $n$.
\end{Proof}

Recall that the categories $\cat\slU$ and $\xc$ satisfy the assumptions of Lemma \ref{two-dual}.
Combining this lemma with Proposition \ref{nemszakad} and Lemma \ref{aut}, we obtain the following.

\begin{Prop}\label{two-dual-spec}
There exists a canonical dual premorphism $\cat{\slU} \to \xc$ if and only if there exists a canonical dual premorphism 
$G^\slU=\cat\slU/G \to \xc/G$.
\end{Prop}

The main results of the section, see Propositions  
\ref{gen2can}, \ref{xc2M}, \ref{nemszakad} and \ref{two-dual-spec},
are summed up in the following theorem.

\begin{Thm}\label{main}
Let $M$ be a quasi-$A$-generated finite-above $E$-unitary inverse monoid 
with $A\subseteq\maxM$, 
put $G=M/\sigma$, and let $\slU$ be a group variety. 
The following statements are equivalent.
\begin{enumerate}
\item \label{main1}
$M$ has an $F$-inverse cover via the group variety $\slU$.
\item \label{main2}
There exists a canonical dual premorphism $G^\slU \to M$.
\item \label{main3} 
There exists a canonical dual premorphism $G^\slU \to\xc/G$.
\item \label{main4} 
There exists a canonical dual premorphism $\cat\slU \to \xc$.
\item  \label{main5}
For any arrow $x$ in $\cat\slU$ and for any $n\in\mathbb{N}_0$, the graph $\Pc{n}(x)$ contains $\tau x$.
\end{enumerate}
\end{Thm}

As an example, we describe a class of non-$F$-inverse finite-above inverse monoids
for which Theorem \ref{main} yields $F$-inverse covers via 
any non-trivial group variety in a straightforward way.
The following observation on the series $\Cc{0}(x),\Cc{1}(x),\ldots$ and $\Pc{0}(x),\Pc{1}(x),\ldots$ of subgraphs plays a crucial role in our argument.
Recall that, given a group variety $\slU$ and a word $w\in \ol{A}^*$, the {\em $\slU$-content $c_\slU(w)$ of $w$} 
consists of the elements $a\in A$ such that $[w]_\slU$ depends on $a$.

\begin{Prop}\label{P0=C0}
\begin{enumerate}
\item \label{C0}
If $x=(\iota x,[p]_\slU,\tau x)$ for some $(\iota x,\tau x)$-path $p$ in $\ol\Gamma$ 
then
$\Cc{0}(x)=\spc{c_\slU(p)}$.
\item \label{xP0=C0}
If $\Cc{0}(x)$ is connected for every arrow $x\in \cat\slU$ then
$\Cc{0}(x) = \Pc{}(x)$ for every $x\in \cat\slU$.
\end{enumerate}
\end{Prop}

\begin{Proof}
The proof of \cite[Lemma 2.1]{szn} can be easily adapted to show (\ref{C0}).
By assumption in (\ref{xP0=C0}), we have 
$\Pc{0}(x)=\Cc{0}(x)$ for any $x\in \cat\slU$.
Applying (\ref{C0}), an easy induction implies that 
$\Cc{n+1}(x)=\Pc{n}(x)$ and $\Pc{n+1}(x)=\Cc{n+1}(x)$ for every $n\in\mathbb{N}_0$ and $x\in \cat\slU$.
This verifies statement (\ref{xP0=C0}).
\end{Proof}

\begin{Ex}
Let $G$ be a group acting on a semilattice $S$  where $S$ has no greatest element, and for every $s \in S$, the set of elements greater than $s$ is finite.
Consider a semidirect product $S \rtimes G$ of $S$ by $G$, and let $M=(S \rtimes G)^1$, the inverse monoid obtained from $S \rtimes G$ by adjoining an identity $1$.
Then $M$ is a finite-above $E$-unitary inverse monoid which is not $F$-inverse, but it has an $F$-inverse cover via any non-trivial group variety.

Notice that $S \rtimes G$ has no identity element, therefore $M\setminus\{1\}=S \rtimes G$.
Recall that the rules of multiplication and taking inverse in $M\setminus\{1\}$ are as follows:
$$(s,g)(t,h)=(s\cdot \acth{g}{t}, gh)\quad \hbox{and}\quad (s,g)^{-1}=(\act{g^{-1}}{s}, g^{-1}).$$
The semilattice of idempotents of $M$ is $(S \times \{1_G\})\cup \{1\}$, and the
natural partial order on $M\setminus\{1\}$  is given by 
$$(s, g) \leq (t,h) \quad \hbox{if and only if} \quad s \leq t \hbox{ and } g=h.$$
The kernel of the projection of $M\setminus\{1\}$ onto $G$, which is clearly a homomorphism, is the least group congruence on $M\setminus\{1\}$.
Hence $M\setminus\{1\}$, and therefore $M$ also is $E$-unitary. 
Morevover, $M$ is finite-above and non-$F$-inverse due to the conditions imposed on $S$.
By Lemma \ref{fa-qAg}, $M$ is quasi-$A$-generated with $A=\maxM$, and it is easy to check that
$\maxM = \max S \times (G\setminus\{1_G\})$ where $\max S$ denotes the maximal elements of $S$.  

Now that all conditions of Theorem \ref{main} are satisfied, construct the graph $\Gamma$: 
its set of vertices is $V_\Gamma=G$ and set of edges is 
\begin{eqnarray*}
E_\Gamma \!\!\!&=&\!\!\!\{(g_1,(s',g),g_2): s' \in \max S\ \hbox{and}\ g_1,g_2,g\in G\ \cr
&&\!\!\!\quad \hbox{such that}\ g\not= 1_G\ \hbox{and}\ g_1g=g_2\},
\end{eqnarray*}
where $\iota (g_1,(s',g),g_2)=g_1$ and $\tau(g_1,(s',g),g_2)=g_2$. 
(This is essentially the Cayley graph of the $A$-generated group $G$ with 
$\epsilon_G\colon A\to G,\ (s',g)\mapsto g$, 
and it is obtained from the Cayley graph of $G$, considered as a $(G\setminus\{1_G\})$-generated group, by replacing each edge with $|\max S|$ copies.)
Let $\slU$ be a non-trivial group variety.
By Proposition \ref{P0=C0}, it suffices to prove that, for each edge $e$ of $\Gamma$, the set of vertices of the graph
$\cl{e}$ is $G$. For, in this case, statement (\ref{C0}) obviously shows that $\Cc{0}(x)$ is connected for every arrow $x$ in $\cat\slU$,
and so statement (\ref{xP0=C0}) implies that Theorem \ref{main}(\ref{main5}) holds for $M$.
Our statement for $M$ follows by the equivalence of Theorem \ref{main}(\ref{main1}) and (\ref{main5}).

Consider an arbitrary edge $e=(g_1,(s',g), g_2) \in E_\Gamma$ and an arbitrary element $h\in G$, and prove that $h$ is a vertex of $\cl{e}$.
Since $g_1$ is obviously a vertex of $\cl{e}$, we can assume that $h\not= g_1$.
Then we have $h=g_1u$ for some $u\in G\setminus\{1_G\}$, and 
$\lambda(e)=(s',g)=(s',u)(s',u)^{-1}(s',g)$.
This implies that $(g_1,(s',u),h)$ is an edge in $\Gamma$ belonging to $\cl{e}$, and so $h$ is, indeed, a vertex of $\cl{e}$.
\end{Ex}

This example sheds light on the generality of 
our construction in contrast with that in \cite{ASzM}. 
By the main result of \cite{szn}, it is known that the 
Margolis--Meakin expansion 
of a group admits an $F$-inverse cover via an Abelian group if and only if the group is cyclic. 
The previous example shows that, 
for any group $G$, there exist finite-above $E$-unitary inverse monoids
with greatest group homomorphic image $G$ that fail to be $F$-inverse but admit $F$-inverse covers via Abelian groups.

\section{$F$-inverse covers via Abelian groups}

In this section, we make further inquiries on how the result of \cite{szn} implying that a Margolis--Meakin expansion of a group 
admits an $F$-inverse cover via an Abelian group if and only
if the group is cyclic generalizes for finite-above $E$-unitary
inverse monoids.
The main result of the section gives a sufficient condition for such an $F$-inverse cover not to exist.

An easy consequence of Theorem \ref{main} is the following:

\begin{Prop}\label{z2}
If $M$ is a finite-above $E$-unitary inverse monoid with $|M/\sigma|\le 2$, then  
$M$ has an $F$-inverse cover via any non-trivial group variety. 
In particular, $M$ has an $F$-inverse cover via an elementary Abelian $p$-group
for any prime $p$.
\end{Prop}

\begin{Proof}
If $|M/\sigma|=1$, that is, $M$ is a semilattice monoid, then $M$ is itself $F$-inverse, 
and the statement holds for any group variety, including the trivial one. 

Now we consider the case $|M/\sigma|=2$. 
Let $A\subseteq \maxM$ such that $M$ is quasi-$A$-generated.
Then the graph $\Gamma$ and the inverse category $\xc$ has two vertices and objects, say, $1$ and $u$.
If $\slU$ is a non-trivial group variety, 
and $q$ is a $(1,u)$-path in $\ol\Gamma$, then $u\not= 1$ implies that $c_\slU(q)$ is non-empty.
Thus $\Cc{0}(x)$ is connected for every arrow $x$ in $\cat{\slU}$, and Proposition \ref{P0=C0}
shows that condition (\ref{main5}) in Theorem \ref{main} is satisfied, completing the proof.
\end{Proof}

This proposition shows that if a finite-above $E$-unitary inverse monoid $M$ has no $F$-inverse cover via an Abelian group
(and consequently, $M$ itself is not $F$-inverse), then $M/\sigma$ has at least two elements distinct from $1$, and
there exists a $\sigma$-class in $M$ containing at least two maximal elements.

From now on, let $M$ be a finite-above $E$-unitary inverse monoid. Let us choose elements $a,b\in M$ with $a\,\sigma\,b$, 
and a $\sigma$-class $v\in M/\sigma$.
Denote by $\max v$ the maximal elements of the $\sigma$-class $v$.
Notice that $\max 1=1_M$, and if $v\not= 1$, then $\max v=v\cap \maxM$.

Consider the following set of idempotents:
$$H(a,b;v)=\{d^{-1}ab^{-1}d:d\in \max v\}.$$
The set of all upper bounds of $H(a,b;v)$ is clearly
$\bigcap\{h^\omega:h\in H(a,b;v)\}$. 
Since $M$ is finite-above, $h^\omega$ is a finite subsemilattice of $E$ for every $h\in H(a,b;v)$ which contains $1_M$.
Therefore $\bigcap\{h^\omega:h\in H(a,b;v)\}$ is also a finite subsemilattice of $E$ containing $1_M$.
This implies that $H(a,b;v)$ has a least upper bound which we denote by $h(a,b;v)$.
The following condition will play a crucial role in this section:

\smallskip

\noindent
{\rm(C)}\qquad
$c\cdot h(a,b;v)\cdot c^{-1}b\not\le a$ for some $c\in \max v$.

\smallskip

Note that if {\rm (C)} is satisfied, then it is not difficult to check that $1,\,u=a\sigma=b\sigma,\ v$
are pairwise distinct elements of $M/\sigma$.
Moreover, $a$ and $b$ are distinct, and $\max v$ contains an element $d$ different from $c$.
Figure 1 shows the arrows of $\slI_M$ related to condition {\rm (C)}.

\begin{figure}
\begin{center}
\includegraphics[width=0.5\linewidth]{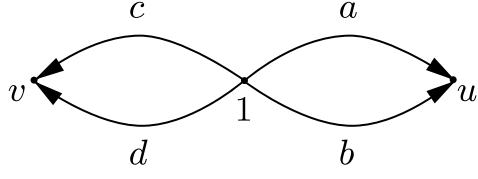}
\end{center}
\caption{The most general constellation of $a,b$, and $v$ in $\slI_M$}
\end{figure}

Denote the variety of Abelian groups by $\slAb$.
The main result of the section is based on the following statement.

\begin{Prop}\label{mprop}
Let $M$ be a finite-above $E$-unitary inverse monoid such that condition {\rm (C)} is satisfied
for some $a,b\in \maxM$ with $a\,\sigma\, b$ and for some $v\in M/\sigma$,
and consider an appropriate $c\in\max v$.
Let $A$ be a quasi-generating set in $M$ such that $A\subseteq \maxM$ and $a,b,c\in A$, 
and consider $M$ as a quasi-$A$-generated inverse monoid.
Then there exists an arrow $x$ in $\cat{\slAb}$ such that $\Pc{1}(x)$ does not contain 
$\tau x$.
\end{Prop}

\begin{Proof}
For every $d\in A$, denote the edge $(1,d,d\sigma)$ of $\Gamma$ by $\ul d$,
and put $u=a\sigma=b\sigma$.
Furthermore, consider the following arrows in $\cat{\slAb}$: 
$$x=(1,[\ul a]_\slAb,u),\quad y=(1,[\ul b]_\slAb,u),\quad z=(1,[\ul c]_\slAb,v).$$
Then we have $z^{-1}xy^{-1}z=(v,[{\ul c}' \ul a {\ul b}' \ul c]_\slAb,v)\in \cat{\slAb}$, where
$[{\ul c}' \ul a {\ul b}' \ul c]_\slAb=[\ul a {\ul b}']_\slAb$ in $F_\slAb(E_\Gamma)$.
For brevity, put $h=h(a,b;v)$, and let
$o$ be a $v$-cycle in $\Gamma^I$ such that 
$\lbl(o)=h$.
It suffices to verify the following two statements:
\begin{equation}\label{c1kovvegl}
\Pc{0}(z^{-1}xy^{-1}z)\subseteq\spc{o},
\end{equation}
\begin{equation}\label{c2kov}
\spc{\ul a}\cap \spc{\ul co{\ul c}'\ul b}\ \hbox{contains no}\ (1,u)\hbox{-path}.
\end{equation}
For, we have $[\ul c({\ul c}'\ul a{\ul b}'\ul c){\ul c}'\ul b]_\slAb=[\ul a]_\slAb$, whence $z(z^{-1}xy^{-1}z)z^{-1}y=x$, and so
$$\Cc1(x)\subseteq \spc{\ul a}\cap(\spc{\ul c}\vee\Pc0(z^{-1}xy^{-1}z)\vee\spc{{\ul c}'\ul b}).$$
Here (\ref{c1kovvegl}) implies
$$\spc{\ul c}\vee\Pc0(z^{-1}xy^{-1}z)\vee\spc{{\ul c}'\ul b}
\subseteq\spc{\ul c}\vee\spc{o}\vee\spc{{\ul c}'\ul b}=\spc{\ul co{\ul c}'\ul b},$$
and so it follows by (\ref{c2kov}) that $\Cc1(x)$ contains no $(1,u)$-path.

Contrary to (\ref{c2kov}), 
assume that the graph $\spc{\ul a}\cap \spc{\ul co{\ul c}'\ul b}$ contains a $(1,u)$-path, say $s$.
Then $\lbl(s)\ge \lbl(\ul a)=a$ and $\lbl(s)\ge \lbl(\ul co{\ul c}'\ul b)=chc^{-1}b$. 
Since $a$ is a maximal element in $M$, the first inequality implies $\lbl(s)=a$, and so
the second contradicts {\rm (C)}.
This shows that (\ref{c2kov}) holds.

To prove 
(\ref{c1kovvegl}), first
we verify that
\begin{equation}\label{c1kov}
\Cc{0}(z^{-1}xy^{-1}z)=\bigcap\bigl\{\spc{t'\ul a{\ul b}'t}:t\ \hbox{is a}\ (1,v)\hbox{-path}\bigr\}.
\end{equation}
It suffices to show that, 
for every $v$-cycle $s$ with $[s]_\slAb=[{\ul c}' \ul a {\ul b}' \ul c]_\slAb=[\ul a {\ul b}']_\slAb$, 
there exists a $(1,v)$-path $t$ such that $\sp{s}=\sp{t'\ul a {\ul b}'t}$.

Let $s$ be a $v$-cycle such that $[s]_\slAb=[\ul a {\ul b}']_\slAb$.
Since $\ul a {\ul b}'$ is a non-trivial simple cycle, 
the former equality implies that $s$ necessarily contains both $\ul a$ and ${\ul b}'$.
Independently of the occurrences of $\ul a$ and $\ul b$ in $s$, the edges $\ul a$ and ${\ul b}'$ appear somewhere
in the $v$-cycle $\tilde s=ss's$ in this order, that is,
$\tilde s=t_0\ul at_1{\ul b}'t_2$ for appropriate paths $t_0,t_1,t_2$.
Moreover, we obviously have
$\sp{\tilde s}=\sp{s}$ and $[\tilde s]_\slAb=[s]_\slAb$.
Putting $\bar s= t_0\ul a{\ul b}'t$, where $t=\ul bt_1{\ul b}'t_2s's$,
we easily see that $\sp{\bar s}=\sp{s}=\sp{t}$ and $[\bar s]_\slAb=[s]_\slAb$.
Finally, the equalities $[t_0\ul a{\ul b}'t]_\slAb=[s]_\slAb=[\ul a{\ul b}']_\slAb$
imply that $[t_0]_\slAb=[t']_\slAb$, and so
$[s]_\slAb=[t'\ul a{\ul b}'t]_\slAb$ and $\sp{s}=\sp{t'\ul a{\ul b}'t}$ follow.
This completes the proof of (\ref{c1kov}).

Turning to the proof of (\ref{c1kovvegl}), assume that 
$k$ is a $v$-cycle in $\Cc{0}(z^{-1}xy^{-1}z)$. 
By (\ref{c1kov}) we see that 
$$\lbl(k)\ge \lbl(t'\ul a{\ul b}'t)=\lbl(t)^{-1}ab^{-1}\lbl(t)$$
for every $(1,v)$-path $t$. 
Since there exists a $(1,v)$-path $t$ with $\lbl(t)=d$ for every $d\in \max v$, 
we obtain that $\lbl(k)$ is an upper bound of $H(a,b;v)$, and so $\lbl(k)\ge h=\lbl(o)$
and $\Pc{0}(z^{-1}xy^{-1}z)\subseteq \Cc{0}(z^{-1}xy^{-1}z) \subseteq \spc{o}$.
This verifies (\ref{c1kovvegl}), and the proof of the proposition is complete. 
\end{Proof}

Combining Proposition \ref{mprop} and Theorem \ref{main}(\ref{main1}) and (\ref{main5}),
we obtain the following sufficient condition for a finite-above $E$-unitary inverse monoid to have no $F$-inverse cover via Abelian groups.

\begin{Thm}
\label{Abelmain}
If $M$ is a finite-above $E$-unitary inverse monoid such that for some $a,b\in \max M$ with $a\,\sigma\, b$ and for some $v\in M/\sigma$,
condition {\rm (C)} is satisfied,
then $M$ has no $F$-inverse cover via Abelian groups.
\end{Thm}


\begin{thebibliography}{99}

\bibitem{A} C.~J.~Ash, {\em Inevitable graphs: A proof of the Type II
conjecture and some related decision procedures},
Internat.~J.~Algebra Comput.~{\bf 1} (1991), 127--146.

\bibitem{ASzM} K.~Auinger and M.~B.~Szendrei, On $F$-inverse covers of inverse monoids,
{\em J.~Pure Appl.~Algebra} {\bf 204} (2006), 493--506.

\bibitem{eil} S.~Eilenberg, {\em Automata, Languages and Machines}, 
Vol. B, Academic Press, New York, 1976.

\bibitem{HR} K.~Henckell and J.~Rhodes, {\em The theorem of Knast,
the ${\bf PG}={\bf BG}$ and type II conjectures}, Monoids and
Semigroups with Applications (Berkeley, CA, 1989), World
Scientific, River Edge, 1991; pp.~453--463.

\bibitem{lo} M.~V.~Lawson, {\em Inverse Semigroups: The Theory of Partial
Symmetries}, World Scientific, Singapore, 1998.

\bibitem{mar-mea} S.~W.~Margolis, J.~C.~Meakin, $E$-unitary inverse monoids and the Cayley graph of a group presentation,
{\em J.~Pure Appl.~Algebra}~{\bf 58} (1989), 45--76.

\bibitem{mar-pin} S.~W.~Margolis, J.-E.~Pin, Inverse semigroups and varieties of finite semigroups, 
{\em J.~Algebra}~{\bf 110} (1987), 306--323.

\bibitem{pet} M.~Petrich, {\em Inverse Semigroups}, Wiley \&\ Sons, New York,
1984. 

\bibitem{szn} N.~Szak\'acs, {\em On the graph condition regarding the F-inverse cover
problem}, to appear in Semigroup Forum, DOI: 10.1007/s00233-015-9713-5

\bibitem{Tilson} B.~Tilson, {\em Categories as algebra: an
essential ingredient in the theory of monoids}, J.~Pure
Appl.~Algebra {\bf 48} (1987), 83--198.

\end{thebibliography}
\end{document}